\documentclass[twoside,11pt,leqno]{article}
\usepackage{amsfonts}

 \def\X{{\cal X}}  \def\H{{\cal H}}

\def\n{{\texttt N}} 
 
\def\r{{\texttt{R}}}

\def\B{B({\cal H})} \def\b{B({\cal X})}

\newtheorem{df}{Definition}[section]
\newtheorem{thm}[df]{Theorem} \newtheorem{pro}[df]{Proposition}
\newtheorem{cor}[df]{Corollary} \newtheorem{ex}[df] {Example}
\newtheorem{rema}[df] {Remark} \newtheorem{lem}[df] {Lemma}
\def\sfstp{{\hskip-1em}{\bf.}{\hskip1em}}

\def\subject#1{\renewcommand{\thefootnote}{}\footnote
{AMS(MOS) subject classification (2010). Primary: {#1}}}

\def\keywords#1{\renewcommand{\thefootnote}{}\footnote
{Keywords: {#1}}}

\def\enddemo{\qed \endtrivlist} \expandafter\let\csname
enddemo*\endcsname=\enddemo

\def\qedsymbol{\ifmmode\bgroup\else$\bgroup\aftergroup$\fi
\vcenter{\hrule\hbox{\vrule
height.5em\kern.5em\vrule}\hrule}\egroup}
\def\qed{\ifmmode\else\unskip\nobreak\fi\quad\qedsymbol}

\title{\bf  Structure of elementary operators defining $m$-left invertible, $m$-selfadjoint and related classes of operators}
\author{\normalsize B.P. Duggal, I.H. Kim}

\date{\centerline{\it \scriptsize Dedicated to Professor Woo Young Lee on the occasion of his 65th birthday}}

\begin{document}

\maketitle \thispagestyle{empty} \vskip-16pt

\subject{Primary47A05, 47A55 Secondary47A80, 47A10.} \keywords{Banach space,  Left/right multiplication operator, elementary operator,  $m$-left invertible, $m$-isometric and $m$-selfadjoint operators, product of operators, perturbation by nilpotents, commuting operators. }
\footnote{The second named author was supported by Basic Science Research Program through the National Research Foundation of Korea (NRF)
 funded by the Ministry of Education (NRF-2019R1F1A1057574).}
\begin{abstract}
 We use elementary, algebraic properties of left, right multiplication operators to prove some deep structural properties of left $m$-invertible, $m$-isometric, $m$-selfadjoint and other related classes of Banach space operators, often adding value to extant results.
\end{abstract}

%%%%%%%%%%%%%%%%%%%%%%%%%%%%%%%%%%%%%%%% SECTION 1

\section {\sfstp Introduction} Let $\b$ (resp., $\B$) denote the algebra of operators on an infinite dimensional complex Banach Space $\X$ (resp., Hilbert space $\H$) into itself.  For $A, B\in\b$, let $L_A$ and $R_B$ $\in B(\b)$ denote, respectively, the operators
$$L_A(X)=AX, \ R_B(X)=XB
$$
of left multiplication by $A$ and right multiplication by $B$. The elementary operators  $\triangle_{A,B}$ and $\delta_{A,B} \in B(\b)$ are then defined by
$$
\triangle_{A,B}(X)=(L_AR_B-I)(X)=L_AR_B(X)-X=AXB-X
$$
and
$$
\delta_{A,B}(X)=(L_A-R_B)(X)=L_A(X)-R_B(X)=AX-XB.
$$
Let $d_{A,B}\in B(\b)$ denote either of the operators $\triangle_{A,B}$ and $\delta_{A,B}$. Let $m\geq 1$ be some integer. The elementary operators
$$
d_{A,B}^m(I)=d_{A,B}(d_{A,B}^{m-1}(I))=\left\{ \begin{array}{l}
\sum_{j=0}^m(-1)^{j}\left(\begin{array}{clcr}m\\j\end{array}\right)A^{m-j}B^{m-j}=0, \ d=\triangle \\ \sum_{j=0}^m(-1)^{j}\left(\begin{array}{clcr}m\\j\end{array}\right)A^{m-j}B^{j}=0, \ d=\delta
\end{array}\right.
$$
have been considered in the recent past by a large number of authors (see \cite{AS1}, \cite{AS2}, \cite{AS3}, \cite{BMN}, \cite{BMN1}, \cite{CKL}, \cite{BD1}, \cite{BD2}, \cite{BD3}, \cite{BM}, \cite{G}, \cite{G1} for further references). Operators satisfying $\triangle_{A,B}^m(I)=0$ have been called left $m$-invertible (\cite{OAM}, \cite{BM}, \cite{G}): $m$-isometric operators $S\in\B$ satisfying $\triangle_{S^*,S}^m(I)=0$ and $(m,C)$-isometric operators  $S\in\B$ satisfying $\triangle_{S^*,CSC}^m(I)=0$ (for some conjugation $C$ of $\H$) are a couple of important example of left $m$-invertible operators. Again, $m$-selfadjoint operators $S\in\B$ satisfying $\delta_{S^*,S}^m(I)=0$ (\cite{Hel}, \cite{TL}) and $(m,C)$-symmetric operators $S\in\B$ satisfying $\delta_{S^*,CSC}^m(I)=0$ (\cite{CKL}) are a couple of important examples of operators $A,\ B$ satisfying $\delta_{A,B}^m(I)=0$.

\medskip

An array of results, amongst them that.
$$
d_{A,B}^m(I)=0 \Longleftrightarrow d_{A,B}^n(I)=0\ {\rm for \ all \ integers}\ n\geq m;
$$
$$
d_{A,B}^m(I)=0 \Longleftrightarrow d_{A^n,B^n}^m(I)=0 \ {\rm for \ all \ integers}\ n\geq 1,
$$
conditions on $A_i, \ B_i \in \b$, $i=1,2$, ensuring
$$
d_{A_1,B_1}^m(I)=0=d_{A_2,B_2}^m(I) \Longrightarrow d_{A_1,A_2,B_1B_2}^{m+n-1}(I)=0;
$$
and conditions on $B$ and an $n$-nilpotent $N$ ensuring
$$
d_{A,B}^m(I)=0 \Longrightarrow d_{A,B+N}^{m+n-1}(I)=0,
$$
is available in extant literature. The authors of these papers have used a wide variety of arguments, amongst them combinatorial arguments, arithematic progressions, (interpolation by) Lagrange polynomials/properties of operator roots of polynomials and the hereditary functional calculus (as developed in \cite{A}).

\bigskip

This paper considers $d_{A,B}^m(I)=0$ from the view point that
$$
d_{A,B}^m(I)=(L_AR_B-1)^m(I)=\left(\sum_{j=0}^m(-1)^{j}\left(\begin{array}{clcr}m\\j\end{array}\right)(L_AR_B)^{m-j}\right)(I)=0, \ d=\triangle ;
$$
$$
d_{A,B}^m(I)=(L_A-R_B)^m(I)=\left(\sum_{j=0}^m(-1)^{j}\left(\begin{array}{clcr}m\\j\end{array}\right)L_A^{m-j}R_B^{j}\right)(I)=0, \ d=\delta
$$
and, using little more them algebraic operations with left-right multiplication operators, provides a uniformly simple argument to prove (often, improved versions of) these results.

\bigskip

The plan of this paper is as follows: Preliminaries are dealt with in Section 2. Here we use elementary arguments to prove some basic properties which operators $A, B$ satisfying $d_{A,B}^m(I)=0$ share with $m$-isometric operators. Amongst other results, we prove here that if $\triangle_{A,B}^m(I)=0$ and $\triangle_{A,B}^{m-1}(I)\neq0$, then
$$
\{ L_A^tR_B^s(I),\ L_A^{t\pm 1}R_B^{s\pm 1}\triangle_{A,B}(I),\ \cdots, \ L_A^{t\pm m-1}R_B^{s\pm m-1}\triangle_{A,B}^{m-1}(I) \}
$$
is a linearly independent set. Sections 3, 4 and 5,  the final three sections (before the section on references), deal with products, perturbations by commuting nilpotents and commuting $A, B$. Here we give elementary proofs for some results which have been proved for $m$-isometric, $(m,C)$-isometric (etc.) operators using a wide variety of often lengthy, occasionally involved, arguments.

\section {\sfstp Preliminaries} Henceforth, unless otherwise stated, $A$ and $B$ shall denote operators in $\b$. We say that the pair $(A,B)\in d^m(I)$ if either $\triangle_{A,B}^{m}(I)=0$ or $\delta_{A,B}^{m}(I)=0$. Since
$$
\triangle_{A,B}^{m+n}(I)=(L_AR_B-I)^n\left(\triangle_{A,B}^m(I)\right)=\sum_{r=0}^n(-1)^{r}\left(\begin{array}{clcr}n\\r\end{array}\right)L_A^{n-r}R_B^{n-r}\left(\triangle_{A,B}^m(I)\right)
$$
and
$$
\delta_{A,B}^{m+n}(I)=(L_A-R_B)^n\left(\delta_{A,B}^m(I)\right)=\sum_{r=0}^n(-1)^{r}\left(\begin{array}{clcr}n\\r\end{array}\right)L_A^{n-r}R_B^{r}\left(\delta_{A,B}^m(I)\right)
$$
for all integers $n\geq1$,
$$
d_{A,B}^m(I)=0 \Longleftrightarrow d_{A,B}^t(I)=0 \ \rm{for \ all \ integers}\ t\geq m.
$$
Again, given an integer $n\geq1$, since
\begin{eqnarray*}
\triangle_{A^n,B^n}^{m}(I)&=&(L_{A^n}R_{B^n}-I)^m(I)=(L_A^nR_B^n-I)^m(I)\\
&=& \left\{(L_{A}R_{B}-I)^m \sum_{j=0}^{m(n-1)}\alpha_j(L_AR_B)^{m(n-1)-j}\right\}(I)\\
&=& \sum_{j=0}^{m(n-1)}\alpha_j(L_AR_B)^{m(n-1)-j}\left(\triangle_{A,B}^m(I)\right)
\end{eqnarray*}
and
\begin{eqnarray*}
\delta_{A^n,B^n}^{m}(I)&=&(L_{A^n}-R_{B^n})^m(I)=(L_A^n-R_B^n)^m(I)\\
&=& \left\{(L_{A}-R_{B})^m \sum_{j=0}^{m(n-1)}\beta_jL_A^{m(n-1)-j}R_B^j\right\}(I)\\
&=& \sum_{j=0}^{m(n-1)}\beta_jL_A^{m(n-1)-j}R^j_B\left(\delta_{A,B}^m(I)\right)
\end{eqnarray*}	
for some scalars $\alpha_j$ and $\beta_j$,
$$
(A,B)\in d^m(I) \Longleftrightarrow (A^n,B^n)\in d^m(I) \ \rm{for \ all \ integers}\ n\geq 1.
$$
(\cite{BMM}, \cite{BM}, \cite{G}).

The hypothesis $(A,B)\in \triangle^m(I)$ implies $B$ is left invertible (and $A$ is right invertible).
Thus, if $(A,B)\in \triangle^m(I)$, $B$ has a dense range and $A$ is injective, then $A, B$ are invertible. Again, if $(A,B)\in \delta^m(I)$ and $B$ is invertible, then
$$
\delta_{A,B}^m(I)=0\Longrightarrow \triangle_{A,B^{-1}}^m(I)=0 \Longrightarrow A \ \rm{is\ right \ invertible}.
$$
Hence, if $B$ is invertible and $A$ is injective, then
$$
\delta_{A,B}^m(I)=0\Longrightarrow A,B \ \rm{ invertible}.
$$
Furthermore, since
$$
(L_{A^{-1}}R_{B^{-1}}-I)^m(I)= (-1)^mL_A^{-m}R_B^{-m}(L_AR_B-I)^m(I)=0,
$$
and
$$
(L_{A^{-1}}-R_{B^{-1}})^m(I)= (-1)^mL_A^{-m}R^{-m}_B(L_A-R_B)^m(I)=0,
$$
we have $(A^{-1},B^{-1})\in d^m(I)$.

\medskip

We say in the following that the pair of operators $(A, B)$ is strictly $d^m(I)$, $(A,B)\in$ strict-$d^m(I)$, if $(A,B)\in d^m(I)$ and $d_{A,B}^{m-1}(I)\neq 0$. The definition implies that if $(A,B)\in$ strict-$d^m(I)$, then $d_{A,B}^r(I)\neq 0$ for all $0\leq r \leq m-1$ and the set
$$
\left\{d_{A,B}^r(I)\right\}_{r=0}^{m-1}=\left\{I, d_{A,B}(I), \cdots, d_{A,B}^{m-1}(I) \right\}
$$
is linearly independent. To see this, assume (to the contrary) that there exist non-zero scalars $a_r$ such that $\sum_{r=0}^{m-1}a_rd_{A,B}^r(I)=0$. Then
$$
0=d_{A,B}^{m-1}\left\{\sum_{r=0}^{m-1}a_rd_{A,B}^r(I)\right\}=a_0d_{A,B}^{m-1}(I) \Longleftrightarrow a_0=0,
$$
$$
0=d_{A,B}^{m-2}\left\{\sum_{r=1}^{m-1}a_rd_{A,B}^r(I)\right\}=a_1d_{A,B}^{m-1}(I) \Longleftrightarrow a_1=0,
$$
$$\cdots$$
$$
0=a_{m-1}d_{A,B}^{m-1}(I) \Longleftrightarrow a_{m-1}=0.
$$

\noindent This being a contradiction, our assertion follows. More is true in the case in which $d=\triangle$.

\begin{lem}\label{lem21} If  $(A,B)\in$ strict-$\triangle^m(I)$, then the sets $\left\{L_A^{t\pm r}\triangle_{A,B}^r(I)\right\}_{r=0}^{m-1}, \left\{R_B^{s\pm r}\triangle_{A,B}^r(I)\right\}_{r=0}^{m-1}$ and $\left\{L_A^{t\pm r}R_B^{s\pm r}\triangle_{A,B}^r(I)\right\}_{r=0}^{m-1}$, where $s,t\geq m-1$ are some integers, are linearly independent.
\end{lem}
\begin{demo} The proof in all cases being similar, we prove the linear independence of the set $\left\{L_A^{t- r}R_B^{s- r}\triangle_{A,B}^r(I)\right\}_{r=0}^{m-1}$.
Suppose to the contrary that there exist scalars $a_r$ such that $\sum_{r=0}^{m-1}a_rL_A^{t- r}R_B^{s- r}\triangle_{A,B}^r(I)=0$. Assume without loss of generality that $s=t+n$ for some integer $n\geq 0$. Then
\begin{eqnarray*}
0&=&\triangle_{A,B}^{m-1}\left(\sum_{r=0}^{m-1}a_rL_A^{t- r}R_B^{s- r}\triangle_{A,B}^r(I)\right) \Longrightarrow a_0A^t\triangle_{A,B}^{m-1}(I)B^{t+n}=0\\
&\Longrightarrow & a_0A^{t+n}\triangle_{A,B}^{m-1}(I)B^{t+n}=0 \Longleftrightarrow a_0\triangle_{A,B}^{m-1}(I)=0  \Longleftrightarrow a_0=0;\\
0&=&\triangle_{A,B}^{m-2}\left(\sum_{r=1}^{m-1}a_rL_A^{t- r}R_B^{s- r}\triangle_{A,B}^r(I)\right) \Longrightarrow a_1A^{t-1}R_B^{s-1}\triangle_{A,B}^{m-1}(I)=0\\
&\Longrightarrow & a_1A^{t+n-1}\triangle_{A,B}^{m-1}(I)B^{t+n-1}=0 \Longleftrightarrow a_1\triangle_{A,B}^{m-1}(I)=0  \Longleftrightarrow a_1=0;\\
& \cdots &\\
0&=&a_{m-1}A^{t-m+1}\triangle_{A,B}^{m-1}(I)B^{s-m+1}\\
& \Longrightarrow& a_{m-1}A^{t+n-m+1}\triangle_{A,B}^{m-1}(I)B^{t+n-m+1}=0\\
&\Longleftrightarrow & a_{m-1}\triangle_{A,B}^{m-1}(I)=0 \Longleftrightarrow a_{m-1}=0.
\end{eqnarray*}	
This is a contradiction.
\end{demo}

A similar result does not hold for $d=\delta$, as the following example shows.

\begin{ex}\label{ex22} Let $\H=\mathbb{C}^3$ (the 3-dimensional complex space) and let $A\in \B$ be the operator
$$
A=\left(\begin{array}{clcr}0 & a & 0\\0 & 0& b\\0& 0& 0\end{array}\right), \ a \ { and} \ b \ { distinct \ non-zero \ reals}.
$$
Let $D$ be the conjugation $D(x,y,z)=({\bar x}, {\bar y}, {\bar z})$.  Then
\begin{eqnarray*}
& &A^3=0,\ DAD=A,\ \delta_{A^*,DAD}^5(I)=0,\ {A^*}^2=\left(\begin{array}{clcr}0 & 0 & 0\\0 & 0& 0\\ab& 0& 0\end{array}\right),\\
& & DA^2D=\left(\begin{array}{clcr}0 & 0 & ab\\0 & 0& 0\\0& 0& 0\end{array}\right)\ {\rm and}\ \delta_{A^*,DAD}^4(I)=\left(\begin{array}{clcr}0 & 0 & 0\\0 & 0& 0\\0& 0& a^2b^2\end{array}\right)\neq 0.
\end{eqnarray*}
Thus $(A^*, DAD) \in$strict-$\delta^5(I)$. However,
$$
L_{A^*}^4\delta_{A^*,DAD}^0(I)= L_{A^*}^3\delta_{A^*,DAD}^1(I)=0,
$$
the set $\left\{L_{A^*}^{4-r}\delta_{A^*,DAD}^r(I)\right\}_{r=0}^4$ is not linearly independent and the set $\left\{\delta_{A^*,DAD}^r(I)\right\}_{r=0}^4$ is  linearly independent.
\end{ex}

We shall require the following lemma in our deliberations.

\begin{lem}\label{lem22} If $R_B^{n-1-i}\triangle_{A,B}^{m+i}(I)=0$ for some integer $n\geq 1$ and $0\leq i \leq n-1$, then $(A,B)\in \triangle^m(I)$.
\end{lem}

\begin{demo}The hypothesis implies $A^{n-1-i}\triangle_{A,B}^{m+i}(I)B^{n-1-i}=0$ for all $0\leq i \leq n-1$. Since $A\triangle_{A,B}^r(I)B=\triangle_{A,B}^{r+1}(I)+\triangle_{A,B}^r(I)$,
$$
A^{n-1-i}\triangle_{A,B}^{m+i}(I)B^{n-1-i}=A^{n-2-i}\triangle_{A,B}^{m+i+1}(I)B^{n-2-i}+A^{n-2-i}\triangle_{A,B}^{m+i}(I)B^{n-2-i};
$$
hence
$$
A^{n-2-i}\triangle_{A,B}^{m+i}(I)B^{n-2-i}=0, \ 0\leq i \leq n-2.
$$
Repeating the argument, another $(n-2)$-times, it follows that $\triangle_{A,B}^m(I)=0$.
\end{demo}

Let ${\cal X}\bar\otimes{\cal X}$ denote the completion, endowed with a reasonable uniform cross-norm, of the algebraic tensor product of $\cal X$ with itself. Let $A\otimes B \in B({\cal X}\bar\otimes{\cal X})$ denote the tensor product operator defined by $A$ and $B$. The following lemma is well known (see \cite{FS}, \cite{HJ} for similar results).

\begin{lem}\label{lem23} If $\sum_{i=1}^nA_i \otimes B_i=0$ for some operators $A_i, B_i \in \b, \ 1\leq i \leq n$ and the set $\left\{B_i \right\}_{1\leq i \leq n}$ is linearly independent, then $A_i=0$ for all $1\leq i\leq n$.
\end{lem}

\section {\sfstp Results: Products}
We start in the following by considering products of commuting pairs of operators satisfying the $d^m(I)$ property. Products of commuting $m$-isometric operators have been considered by Berm\'{u}dez {\em et al} (\cite{BMN} and   \cite{TL}), of commuting left $m$-invertible operators by \cite{G}, of tensor products of left $m$-invertible operators by Duggal and Muller (\cite{BM}) and of commuting $m$-selfadjoint operators by Trieu Le (\cite{TL}). The following theorem employs an elementary algebraic argument to prove these results for operators satisfying the $d^m(I)$ property. Let $[A,B]=AB-BA$ denote the commutator of $A$ and $B$.

\begin{thm}\label{thm23} Let $A_i, B_i \in \b, i=1,2,$ be such that $[A_1, A_2]=0=[B_1,B_2]$.

\vskip4pt\noindent (i) If $(A_1,B_1)\in d^m(I)$ and $(A_2,B_2)\in d^n(I)$, then $(A_1A_2, B_1B_2)\in d^{m+n-1}(I)$.

\vskip4pt\noindent (ii) If $(A_1,B_1)\in \triangle^m(I)$ and $(A_2,B_2)\in \triangle^n(I)$, then $(A_1A_2, B_1B_2)\in$ strict-$\triangle^{m+n-1}(I)$ if and only if $(L_{A_2}R_{B_2})^{m-1}\triangle_{A_1,B_1}^{m-1}(I)\triangle_{A_2,B_2}^{n-1}(I)\neq 0$, equivalently, $\triangle_{A_1,B_1}^{m-1}(I)\triangle_{A_2,B_2}^{n-1}(I)\neq 0$.

\vskip4pt\noindent (iii) If $(A_1,B_1)\in \delta^m(I)$ and $(A_2,B_2)\in \delta^n(I)$, then $(A_1A_2, B_1B_2)\in$ strict-$\delta^{m+n-1}(I)$ if and only if $R_{B_1}^{n-1}L_{A_2}^{m-1}\delta_{A_1,B_1}^{m-1}(I)\delta_{A_2,B_2}^{n-1}(I)=L_{A_1}^{n-1}R_{B_2}^{m-1}\delta_{A_1,B_1}^{m-1}(I)\delta_{A_2,B_2}^{n-1}(I)\neq 0$.
\end{thm}

\begin{demo} (i) The hypothesis $[A_1, A_2]=0=[B_1,B_2]$ implies $L_{A_1}, L_{A_2}, R_{B_1}$ and $R_{B_2}$ all (mutually) commute. Since
\begin{eqnarray*}
(L_{A_1} L_{A_2} R_{B_1}R_{B_2}-I)^r&=&\left\{L_{A_2} (L_{A_1} R_{B_1}-I)R_{B_2}+(L_{A_2} R_{B_2}-I)\right\}^r\\
&=& \sum_{j=0}^r\left(\begin{array}{clcr}r\\j\end{array}\right)(L_{A_2}R_{B_2})^{r-j}(L_{A_1}R_{B_1}-I)^{r-j}(L_{A_2}R_{B_2}-I)^{j}
\end{eqnarray*}
and
$$
\left\{(L_{A_1} R_{B_1}-I)^{r-j}(L_{A_2} R_{B_2}-I)^j\right\}(I)=\triangle_{A_2,B_2}^{j}(I)\triangle_{A_1,B_1}^{r-j}(I)=\triangle_{A_1,B_1}^{r-j}(I)\triangle_{A_2,B_2}^{j}(I)
$$
for all integers $j\geq 1$, we have
$$
\triangle_{A_1A_2,B_1B_2}^{r}(I)=\sum_{j=0}^r\left(\begin{array}{clcr}r\\j\end{array}\right)(L_{A_2}R_{B_2})^{r-j}\triangle_{A_2,B_2}^{j}(I)\triangle_{A_1,B_1}^{r-j}(I).
$$
Let $r=m+n-1$. Then $\triangle_{A_2,B_2}^{j}(I)=0$ for all $j\geq n$, and since
$$
j\leq n-1 \Longrightarrow m+n-1-j\geq m \Longrightarrow \triangle_{A_1,B_1}^{m+n-1-j}(I)=0,
$$
we have
$$
(A_1A_2, B_1B_2)\in \triangle^{m+n-1}(I).
$$
Considering now
\begin{eqnarray*}
(L_{A_1} L_{A_2} -R_{B_1}R_{B_2})^r&=&\left\{L_{A_2} (L_{A_1}- R_{B_1})+(L_{A_2}- R_{B_2})R_{B_1}\right\}^r\\
&=& \sum_{j=0}^r\left(\begin{array}{clcr}r\\j\end{array}\right)L_{A_2}^{r-j}R_{B_1}^j(L_{A_1}-R_{B_1})^{r-j}(L_{A_2}-R_{B_2})^{j},
\end{eqnarray*}
we have
$$
\delta_{A_1A_2,B_1B_2}^r(I)=\sum_{j=0}^r\left(\begin{array}{clcr}r\\j\end{array}\right)L_{A_2}^{r-j}R_{B_1}^j\delta_{A_1,B_1}^{r-j}(I)\delta_{A_2,B_2}^{j}(I).
$$
Hence, letting $r=m+n-1$ and recalling $\delta_{A_1,B_1}^{m}(I)=0=\delta_{A_2,B_2}^{n}(I)$, we have
$$
\delta_{A_1A_2,B_1B_2}^{m+n-1}(I)=\sum_{j=0}^{m+n-1}\left(\begin{array}{clcr}m+n-1\\j\end{array}\right)R_{B_1}^j\delta_{A_1,B_1}^{m+n-1-j}(I)L_{A_2}^{m+n-1-j}\delta_{A_2,B_2}^{j}(I)=0.
$$

(ii) If $\triangle_{A_1,B_1}^{m}(I)=0=\triangle_{A_2,B_2}^{n}(I)$, then
\begin{eqnarray*}
\triangle_{A_1A_2,B_1B_2}^{m+n-2}(I)&=&\sum_{j=0}^{m+n-2}\left(\begin{array}{clcr}m+n-2\\j\end{array}\right)(L_{A_2}R_{B_2})^{m+n-2-j}
\triangle_{A_2,B_2}^{j}(I)\triangle_{A_1,B_1}^{m+n-2-j}(I)\\
&=&\left(\begin{array}{clcr}m+n-2\\m-1\end{array}\right)(L_{A_2}R_{B_2})^{m-1}\triangle_{A_2,B_2}^{n-1}(I)\triangle_{A_1,B_1}^{m-1}(I)\\
&=&\left(\begin{array}{clcr}m+n-2\\m-1\end{array}\right)\triangle_{A_2,B_2}^{n-1}(I)\triangle_{A_1,B_1}^{m-1}(I)
\end{eqnarray*}
(since $\triangle_{A_2,B_2}^{n}(I)=0$ implies $A_2^{m-1}\triangle_{A_2,B_2}^{m-1}(I)B_2^{m-1}=\triangle_{A_2,B_2}^{n-1}(I))$. Hence
\begin{eqnarray*}
\triangle_{A_1A_2,B_1B_2}^{m+n-2}(I)\neq 0 &\Longleftrightarrow & (L_{A_2}R_{B_2})^{m-1}\triangle_{A_2,B_2}^{n-1}(I)\triangle_{A_2,B_2}^{m-1}(I)\neq 0\\
&\Longleftrightarrow&\triangle_{A_2,B_2}^{n-1}(I)\triangle_{A_1,B_1}^{m-1}(I)\neq 0
\end{eqnarray*}
(equivalently, $(A_1, B_1)\in$ strict-$\triangle^m(I)$ and $(A_2, B_2)\in$ strict-$\triangle^n(I))$.

(iii) Again, since $\delta_{A_1,B_1}^{m}(I)=0=\delta_{A_2,B_2}^{n}(I),$
\begin{eqnarray*}
\delta_{A_1A_2,B_1B_2}^{m+n-2}(I)&=& \sum_{j=0}^{m+n-2}\left(\begin{array}{clcr}m+n-2\\j\end{array}\right)L_{A_2}^{m+n-2-j}R_{B_1}^j\delta_{A_1,B_1}^{m+n-2-j}(I)\delta_{A_2,B_2}^{j}(I)\\
&=&\left(\begin{array}{clcr}m+n-2\\n-1\end{array}\right)L_{A_2}^{m-1}R_{B_1}^{n-1}\delta_{A_1,B_1}^{m-1}(I)\delta_{A_2,B_2}^{n-1}(I).
\end{eqnarray*}
Evidently,
\begin{eqnarray*}
\delta_{A_2,B_2}^{n}(I)=0 &\Longleftrightarrow &L_{A_2}\delta_{A_2,B_2}^{n-1}(I)=R_{B_2}\delta_{A_2,B_2}^{n-1}(I)\\
&\Longrightarrow & L_{A_2}^{m-1}\delta_{A_2,B_2}^{n-1}(I)=R_{B_2}^{m-1}\delta_{A_2,B_2}^{n-1}(I)
\end{eqnarray*}
and
\begin{eqnarray*}
\delta_{A_1,B_1}^{m}(I)=0 &\Longleftrightarrow & L_{A_1}\delta_{A_1,B_1}^{m-1}(I)=R_{B_1}\delta_{A_1,B_1}^{m-1}(I)\\
&\Longrightarrow & R_{B_1}^{n-1}\delta_{A_1,B_1}^{m-1}(I)=L_{A_1}^{n-1}\delta_{A_1,B_1}^{m-1}(I).
\end{eqnarray*}
Hence
\begin{eqnarray*}
(A_1A_2, B_1B_2)\in \ {\rm strict-}\delta^{m+n-1}(I)  &\Longleftrightarrow & L_{A_2}^{m-1}R_{B_1}^{n-1}\delta_{A_1,B_1}^{m-1}(I)\delta_{A_2,B_2}^{n-1}(I)\neq 0\\
&\Longleftrightarrow&L_{A_1}^{n-1}R_{B_2}^{m-1}\delta_{A_1,B_1}^{m-1}(I)\delta_{A_2,B_2}^{n-1}(I)\neq 0.
\end{eqnarray*}
\end{demo}

The condition that $\delta_{A_1,B_1}^{m-1}(I)\neq 0 \neq \delta_{A_2,B_2}^{n-1}(I)$, although necessary for $(A_1A_2, B_1B_2)\in$ strict-$\delta^{m+n-1}(I)$ in Theorem \ref{thm23}(iii), is not sufficient.

\begin{ex}\label{ex24} Let $A\in \B, \ \H={\mathbb C}^3$, be the operator of Example \ref{ex22} and let the conjugation  $C$  of $\H$ be defined by $C(x,y,z)=(\bar{z},\bar{y},\bar{x})$. Then $[C,D]=0, \ (A^*, CAC)\in$strict-$\delta^3(I)$, $(A^*, DAD)\in$strict-$\delta^5(I)$, $({A^*}^2, CACDAD)\in$strict-$\delta^7(I)$ and $${A^*}^2D{A^*}^4D\delta_{A^*, CAC}^2(I)\delta_{A^*,DAD}^4(I)=0$$ (even though neither of the operators $\delta_{A^*, CAC}^2(I)$ and $\delta_{A^*,DAD}^4(I)$ is $0$).
\end{ex}

Translated to tensor products, Theorem \ref{thm23} implies:

\begin{pro}\label{pro25} Let $A_i, B_i \in \b, i=1,2.$

\vskip4pt\noindent (i) If $(A_1,B_1)\in d^m(I)$ and $(A_2,B_2)\in d^n(I)$, $d=\triangle \ \rm{or} \ \delta$ (exclusive 'or' as always), then $(A_1\otimes A_2. B_1\otimes B_2)\in d^{m+n-1}(I\otimes I)$.

\vskip4pt\noindent (ii) Any two of the hypotheses (a) $(A_1, B_1)\in {\rm strict}-\triangle^m(I)$; (b) $(A_2, B_2)\in {\rm strict}-\triangle^n(I)$; (c) $(A_1\otimes A_2, B_1\otimes B_2)\in {\rm strict}-\triangle^{m+n-1}(I\otimes I)$ implies the other.

\vskip4pt\noindent (iii) If $d=\delta$ in part (i), then $(A_1\otimes A_2. B_1\otimes B_2)\in {\rm strict}-\triangle^{m+n-1}(I\otimes I)$ if and only if
$$
(A_1^{n-1}\otimes A_2^{m-1})\delta_{A_1\otimes I, I\otimes B_1}^{m-1}(I\otimes I)\delta_{A_2\otimes I, I\otimes B_2}^{n-1}(I\otimes I)\neq 0.
$$
\end{pro}

 \begin{demo} Define $S_i, T_i \in B({\cal X}\otimes {\cal X}), \ i=1,2$, by $S_1=(A_1 \otimes I),\ S_2=(I \otimes A_2), \ T_1=(B_1 \otimes I)$ and $ T_2=(I \otimes B_2)$. Then $[S_1, S_2]=0=[T_1,T_2], \ S_1S_2=(A_1 \otimes A_2)$ and $T_1T_2=(B_1 \otimes B_2)$. Theorem \ref{thm23} applies and we have:

\vskip3pt\noindent (i) $d_{S_1S_2, T_1T_2}^{m+n-1}(I)=0$ implies $d_{A_1\otimes A_2, B_1\otimes B_2}^{m+n-1}(I \otimes I)=0$.

\vskip3pt\noindent  (ii) Since $(A_1\otimes A_2, B_1\otimes B_2)\in$ strict-$\triangle^{m+n-1}(I \otimes I)$ if and only if $\triangle_{A_1 \otimes I, B_1 \otimes I}^{m-1}(I \otimes I)\neq 0 \neq \triangle_{I \otimes A_2, I \otimes B_2}^{n-1}(I \otimes I)$, $(A_1, B_1)\in$strict-$\triangle^m(I)$ if and only if $\triangle_{A_1,B_1}^{m-1}(I)\neq 0$ (equivalently,  $\triangle_{A_1 \otimes I,B_1 \otimes I}^{m-1}(I \otimes I)\neq 0$) and $(A_2, B_2)\in$strict-$\triangle^n(I)$ if and only if $\triangle_{A_2,B_2}^{m-1}(I)\neq 0$ (equivalently,  $\triangle_{I \otimes A_2, I \otimes B_2}^{n-1}(I \otimes I)\neq 0$), any two of (a), (b) and (c) imply the third.

\vskip3pt\noindent The proof of (iii) being evident, the proof is complete.
 \end{demo}

\section {\sfstp Results: Perturbation by commuting nilpotents} Perturbation by commuting nilpotents of $m$-isometric, $m$-selfadjoint and left $m$-invertible (etc.) operators has been considered by a number of authors in the recent past, amongst them Berm\'{u}dez {\em et al} (\cite{BMN1}, \cite{BMMN}, Le \cite{TL} and  Gu \cite{G}). The following theorem provides a particularly elementary proof of the currently available  results for pairs $(A,B)\in d^m(I)$.

\begin{thm}\label{thm31} If $(A,B)\in d^m(I)$, and $N$ is an $n$-nilpotent which commutes with $B$, then $(A,B+N)\in d^{m+n-1}(I)$. Furthermore, $(A,B+N)\in$ strict-$ d^{m+n-1}(I)$ if and only if
$$
R_N^{n-1}\delta_{A,B}^{m-1}(I)\neq 0, \ d=\delta
$$
and
$$
L_A^{n-1}R_N^{n-1}\triangle_{A,B}^{m-1}(I)\neq 0, \ d=\triangle.
$$
\end{thm}

\begin{demo} A straightforward calculation shows that
\begin{eqnarray*}
\triangle_{A,B+N}^{m+n-1}(I)&=& \left\{(L_AR_B-I)+L_AR_N \right\}^{m+n-1}(I)\\
&=&\left\{ \sum_{j=0}^{m+n-1}\left(\begin{array}{clcr}m+n-1\\j\end{array}\right)(L_{A}R_N)^j(L_AR_B-I)^{m+n-1-j}\right\}(I)\\
&=&\sum_{j=0}^{m+n-1}\left(\begin{array}{clcr}m+n-1\\j\end{array}\right)(L_{A}R_N)^j\triangle_{A,B}^{m+n-1-j}(I)
\end{eqnarray*}
and
\begin{eqnarray*}
\delta_{A,B+N}^{m+n-1}(I)&=& \left\{(L_A-R_B)-R_N \right\}^{m+n-1}(I)\\
&=&\left\{ \sum_{j=0}^{m+n-1}(-1)^j\left(\begin{array}{clcr}m+n-1\\j\end{array}\right)R_N^j(L_A-R_B)^{m+n-1-j}\right\}(I)\\
&=&\sum_{j=0}^{m+n-1}(-1)^j\left(\begin{array}{clcr}m+n-1\\j\end{array}\right)R_N^j\delta_{A,B}^{m+n-1-j}(I).
\end{eqnarray*}

The operator $N$ being $n$-nilpotent, $N^j=0$ for all $j\geq n$. For $j\leq n-1$ (equivalently, $-j\geq -n+1$), $m+n-1-j \geq m$. Hence $(A,B+N)\in d^{m+n-1}(I)$.

Since $(A,B+N)\in$ strict-$d^{m+n-1}(I)$ if and only if $d_{A,B+N}^{m+n-2}(I)\neq 0$, it follows that  $(A,B+N)\in$ strict-$d^{m+n-1}(I)$ if and only if
$$
\left(\begin{array}{clcr}m+n-1\\n-1\end{array}\right)(L_{A}R_N)^{n-1}\triangle_{A,B}^{m-1}(I)\neq 0, \ d=\triangle
$$
and
$$
\left(\begin{array}{clcr}m+n-1\\n-1\end{array}\right)R_N^{n-1}\delta_{A,B}^{m-1}(I)\neq 0, \ d=\delta.
$$
The proof is now evident.
\end{demo}

Evidently, $(A,B)\in d^m(I)$ implies  $(B^*,A^*)\in d^m(I^*)$, and if $M$ is an $n_1$-nilpotent which commutes with $A$, then $(B^*,A^*+M^*)\in d^{m+n_1-1}(I^*)$. Hence, Theorem \ref{thm31} implies:

\begin{cor}\label{cor32} If $(A,B)\in d^m(I)$, and $N_i \ (i=1,2)$ are $n_i$-nilpotents such that $[A,N_2]=0=[B,N_1]$, then $(A+N_2,B+N_1)\in d^{m+n_1+n_2-2}(I).$
\end{cor}

If $A, B$ and $N_i$ ($i=1,2$) are the operators of Corollary \ref{cor32}, then
$$
\triangle_{A+N_2,B+N_1}^r(I)=\left\{ \sum_{j=0}^{r}\left(\begin{array}{clcr}r\\j\end{array}\right)\triangle_{A,B+N_1}^{r-j}(L_{N_2}R_{B+N_1})^j\right\}(I)
$$
 and
$$
\delta_{A+N_2,B+N_1}^r(I) =\left\{ \sum_{j=0}^{r}\left(\begin{array}{clcr}r\\j\end{array}\right)\delta_{A,B+N_1}^{r-j}L_{N_2}^j\right\}(I)
$$
and it follows from an argument similar to the one above that $(A+N_2,B+N_1)\in$ strict-$ d^{m+n_1+n_2-2}(I)$ if and only if
$$
\triangle_{A,B+N_1}^{m+n_1-2}(I)(L_{N_2}R_{B+N_1})^{n_2-1}\neq 0 , \ d=\triangle
$$
and
$$
\delta^{m+n_1-2}_{A,B+N_1}(I)L^{n_2-1}_{N_2}\neq 0, \ d=\delta.
$$

In the particular case in which $m=1$,  Theorem \ref{thm31} implies that a necessary and sufficient condition for $(A,B+N)\in$ strict-$\triangle^n(I)$ is that $A^{n-1}N^{n-1}B^{n-1}=A^{n-1}B^{n-1}N^{n-1}=N^{n-1}\neq 0$. Hence :

 \begin{cor}\label{cor33} If $(A,B)\in \triangle(I)$, and $N$ is an $n$-nilpotent which commutes with $B$, then $(A,B+N)\in$ strict-$\triangle^n(I)$.
\end{cor}

An $m$-isometric version of Corollary \ref{cor33} was first proved by Berm\'{u}dez {\em et.al.} (\cite{BMN1}, Theorem 2.2). Observe that if $d=\delta$ and $m=1$, then $\delta_{A,B}(I)=0$ if and only if $A=B$ and then $(A,B+N)\in$ strict-$\delta^n(I)$ in Theorem \ref{thm31}.

If $m>1$, then the hypothesis  $(A,B)\in$ strict-$d^m(I)$ is not enough to guarantee  $(A,B+N)\in$ strict-$d^{m+n-1}(I)$, as the following (slightly changed) example from \cite[Example 7]{G} shows.

\begin{ex}\label{ex34} Define $A,B \in B({\cal X} \oplus {\cal X})$ by $A=I \oplus A_1$ and $B=I \oplus B_1$, where $(A_1,B_1)\in$ strict-$d^m(I)$. Let $N_1 \in \b$ be an $n$-nilpotent operator, and define $N \in B({\cal X} \oplus {\cal X})$ by $N=N_1 \oplus 0$. Then neither of the operators $d_{A,B}^{m-1}(I ),\ N^{n-1}$ and $A^{n-1}$ is  the $0$ operator,
\begin{eqnarray*}
& &d^r_{A, B+N}(I\oplus I)=d^r_{I, I+N_1}(I)\oplus d^r_{A_1, B_1}(I),\\
 & &(I, I+N_1)\in {\rm strict}-d^n(I), \ (A_1 ,B_1)\in {\rm strict}-\triangle^m(I), \ {\rm and}\\
& & (A, B+N) \in {\rm strict}-d^t(I), \ t={\rm max}\{n, m\}< m+n-1.
\end{eqnarray*}
Observe that
\begin{eqnarray*} L_A^{n-1}R_N^{n-1}\triangle^{m-1}_{A,B}(I\oplus I)&=&(I\oplus A^{n-1}_1)(0\oplus \triangle^{m-1}_{A_1,B_1}(I))(N^{n-1}_1\oplus 0)\\ &=& 0\\
&=&(0\oplus\delta^{m-1}_{A_1,B_1}(I))(N^{n-1}_1\oplus 0)=R^{n-1}_N\delta^{m-1}_{A,B}(I\oplus I)
\end{eqnarray*}
({\em cf.} Theorem \ref{thm31}).
 \end{ex}
Clearly, $t={\rm max}\{n, m\}=n$ in (the example above in) the case in which $m=1$. The case $d=\delta$ of Theorem \ref{thm31} is a particular case of the following proposition.

\begin{pro}\label{pro35} If $A_i, B_i \in \b$ are such that $[A_1, A_2]=0=[B_1, B_2], \ (A_1, B_1)\in \delta^m(I)$ and $(A_2, B_2)\in \delta^n(I)$, then $(A_1+A_2, B_1+B_2)\in$ strict-$ \delta^{m+n-1}(I)$ if and only if  $\delta_{A_1,B_1}^{m-1}(I)\delta_{A_2,B_2}^{n-1}(I)\neq 0$.
\end{pro}

\begin{demo} If we let $S=(L_{A_1}+L_{A_2}-R_{B_1}-R_{B_2})(I)=\left\{(L_{A_1}-R_{B_1})+(L_{A_2}-R_{B_2})\right\}(I)$, then
\begin{eqnarray*}
S^{m+n-1}&=&\sum_{j=0}^{m+n-1}\left(\begin{array}{clcr}m+n-1\\j\end{array}\right)\delta_{A_1,B_1}^{m+n-1-j}(I)\delta_{A_2,B_2}^{j}(I) \\ &=& \sum_{j=0}^{m+n-1}\left(\begin{array}{clcr}m+n-1\\j\end{array}\right)\delta_{A_2,B_2}^{j}(I)\delta_{A_1,B_1}^{m+n-1-j}(I).
\end{eqnarray*}
The hypothesis $\delta_{A_1,B_1}^{m}(I)=0=\delta_{A_2,B_2}^{n}(I)$ forces $S^{m+n-1}=0$. Since
\begin{eqnarray*}
S^{m+n-2}&=&\sum_{j=0}^{m+n-2}\left(\begin{array}{clcr}m+n-2\\j\end{array}\right)\delta_{A_1,B_1}^{m+n-2-j}(I)\delta_{A_2,B_2}^{j}(I)\\
&=& \left(\begin{array}{clcr}m+n-2\\n-1\end{array}\right)\delta_{A_1,B_1}^{m-1}(I)\delta_{A_2,B_2}^{n-1}(I),
\end{eqnarray*}
it follows that $S^{m+n-2}\neq 0$ if and only if $\delta_{A_1,B_1}^{m-1}(I)\delta_{A_2,B_2}^{n-1}(I) \neq  0$.
\end{demo}

The case $d=\delta$ of Theorem \ref{thm31} is obtained from Proposition \ref{pro35} upon choosing $A_2=0$ and $B_2=N$ (for then $\delta_{0,B_2}^{n}(I)=0$ if and only if $B_2^n=0$). It is clear from Example \ref{ex34} that the condition $\delta^{m-1}_{A_,B_1}(I)\neq 0\neq \delta^{n-1}_{A_2,B_2}(I)$ in Proposition \ref{pro35} though necessary is in no way sufficient for $(A_1+A_2,B_1+B_2)\in$ strict-$\delta^{m+n-1}(I)$. If we let $B_i=S_i^*$ and $A_i=S_i, \ i=1,2$, for some operators $S_i\in \B$, and choose $n=1$, then Proposition \ref{pro35} says: $S_1+S_2$ is $m$-selfadjoint (\cite{TL}, Corollary 2.9). Furthermore, $S_1+S_2\in$ strict-$m$-selfadjoint if and only if $S_1\in$ strict-$m$-selfadjoint.

\begin{rema}\label{rem36} $A\in\B$ is said to be an $n$-Jordan operator (or, a Jordan operator of order $n$) if $A=S+N$ for some selfadjoint operator $S\in\B$ (i.e., for some $S\in\B$ satisfying $\delta_{S^*,S}(I)=0$) and an $n$-nilpotent operator $N\in\B$ such that $[S,N]=0$. Theorem \ref{thm31} says that "if $\delta_{S^*,S}(I)=0$, $[S,N]=0$, $N^n=0$ and $A=S+N$, then $A$ is $(2n-1)$-Jordan and satisfies $\delta^{2n-1}_{A^*,A}(I)=0$" (\cite{McR}, Theorem 3.2).
\end{rema}

\begin{rema}\label{rem37} Bayart (\cite{FB}) defines a Banach space operator $S\in \b$ to be $m$-isometric if $\sum_{j=0}^{m}(-1)^j\left(\begin{array}{clcr}m\\j\end{array}\right)\|S^{m-j}x\|^2=0$ for all $x\in {\cal X}$. Hilbert space $m$-isometric operators are left $m$-invertible operators, but a similar description does not fit the (generalized) Banach space definition. Arguments similar to those used to prove Theorems \ref{thm23} and \ref{thm31} do not extend to prove similar results for Banach space $m$-isometries. This in a way is to be expected: Whereas it is true that if $A, B \in\b$ are commuting operators such that $A$ is $m$-isometric and $B$ is $n$-isometric  then $AB$ is $(m+n-1)$-isometric (\cite{BMN}, Theorem 3.3), the perturbation by a commuting $n$-nilpotent operator of an $m$-isometric (Banach space) operator may not be an $(m+2n-2)$-isometric operator (\cite{BMMN}, Example 4.1). Observe that Bayart's definition of $m$-isometries converts the problem of commuting products into a problem in arithmetic progressions.
\end{rema}

Let $N\in \b$ be an $n$-nilpotent operator, and let $S^r, T^r \in B({\cal X} \bar{\otimes} {\cal X}), \ r\geq 1$ some integer, be the operators
 \begin{eqnarray*}
S^r&=&\triangle_{A \otimes I +I \otimes N, B \otimes I}^r(I \otimes I) =\left\{(L_{A \otimes I}R_{B \otimes I}-I \otimes I)+L_{A \otimes I}R_{I \otimes N}\right\}^r(I \otimes I)\\
&=&\left\{\sum_{j=0}^{r}\left(\begin{array}{clcr}r\\j\end{array}\right)(L_{I \otimes N}R_{B \otimes I})^j(L_{A \otimes I}R_{B \otimes I}-I \otimes I)^{r-j}\right\}(I \otimes I)\\
&=&\sum_{j=0}^{r}\left(\begin{array}{clcr}r\\j\end{array}\right)\triangle_{A,B}^{r-j}(I)B^j\otimes N^j
\end{eqnarray*}
and
\begin{eqnarray*}
T^r&=&\delta_{A \otimes I +I \otimes N, B \otimes I}^r(I \otimes I) =\left\{\sum_{j=0}^{r}\left(\begin{array}{clcr}r\\j\end{array}\right)(L_{A \otimes I}-R_{B \otimes I})^{r-j}L_{I \otimes N}^j\right\}(I \otimes I)\\
&=&\sum_{j=0}^{r}\left(\begin{array}{clcr}r\\j\end{array}\right) \delta_{A,B}^{r-j}(I) \otimes N^j.
\end{eqnarray*}

\noindent  Let $E=A \otimes I, \ F=B \otimes I$ and $Q=I \otimes N$. Then $[E,Q]=0=Q^n$. If we now assume that $(A,B)\in$strict-$d^m(I)$, then it follows from an application of Theorem \ref{thm31} that
 $$
 (E+Q, F)\in d^{m+n-1}(I \otimes I),
$$
$$
 \triangle_{E+Q, F}^{m+n-2}(I \otimes I)=S^{m+n-2}\neq 0 \neq T^{m+n-2}=\delta_{E+Q, F}^{m+n-2}(I \otimes I)
 $$
 and
\begin{eqnarray} \label{eqnarray1} \triangle_{A,B}^{m-1}(I)B^{n-1} \otimes N^{n-1} \neq 0 \neq \delta_{A,B}^{m-1}\otimes N^{n-1}.\end{eqnarray}

The following proposition is proved in (\cite{G}, Theorem 18)  for the case in which $d=\triangle$; our proof, however, differs from that in (\cite{G}).

\begin{pro}\label{pro37} Given operators $A, B, N \in \b$, any two of the following conditions implies the third.

\vskip4pt\noindent (i) $(A,B)\in$ strict-$d^m(I)$.

\vskip4pt\noindent (ii) $N^n=0$.

\vskip4pt\noindent (iii) $(A \otimes I +I \otimes N, B \otimes I)\in$ strict-$d^{m+n-1}(I \otimes I)$.
\end{pro}

\begin{demo} $(i) \land (ii) \Longrightarrow (iii)$: We have already seen (above) that $(A \otimes I +I \otimes N, B \otimes I)\in d^{m+n-1}(I \otimes I)$. Since $\triangle_{A,B}^{m-1}(I)B^{n-1}=0$ implies $A^{n-1}\triangle_{A,B}^{m-1}(I)B^{n-1}=\triangle_{A,B}^{m-1}(I)=0$ (Recall: $\triangle_{A,B}^{m}(I)=A\triangle_{A,B}^{m-1}(I)B-\triangle_{A,B}^{m-1}(I)=0)$, and since $N^{n-1}\neq 0\neq d_{A,B}^{m-1}(I)$, $(1)$ holds.

 \bigskip

\noindent $(i) \land (iii) \Longrightarrow (ii)$: Hypothesis (i) implies
$$
S^{m+n-1}=\sum_{j=n}^{m+n-1}\left(\begin{array}{clcr}m+n-1\\j\end{array}\right)\triangle_{A,B}^{m+n-1-j}(I)B^j\otimes N^j
$$
and
$$
T^{m+n-1}=\sum_{j=n}^{m+n-1}\left(\begin{array}{clcr}m+n-1\\j\end{array}\right)\delta_{A,B}^{m+n-1-j}(I)\otimes N^j.
 $$
 Here the sets $\left\{\triangle_{A,B}^{m+n-1-j}(I)B^j \right\}_{j=n}^{m+n-1}$ and $\left\{\delta_{A,B}^{m+n-1-j}(I) \right\}_{j=n}^{m+n-1}$ are linearly independent. Hence we must have $N^j=0$ for all $n\leq j\leq m+n-1$.

\bigskip

\noindent $(ii) \land (iii) \Longrightarrow (i)$: Hypotheses $(ii) \land (iii)$ imply
$$
S^{m+n-1}=\sum_{j=0}^{n-1}\left(\begin{array}{clcr}m+n-1\\j\end{array}\right)\triangle_{A,B}^{m+n-1-j}(I)B^j \otimes N^j=0
$$
and
$$
T^{m+n-1}=\sum_{j=0}^{n-1}\left(\begin{array}{clcr}m+n-1\\j\end{array}\right)\delta_{A,B}^{m+n-1-j}(I)\otimes N^j=0.
$$
The set $\left\{N^j \right\}_{j=0}^{n-1}$ being linearly independent, we must have $\triangle_{A,B}^{m+j}(I)B^j=0=\delta_{A,B}^{m+j}(I)$ for all $0\leq j\leq n-1$. Applying Lemma 2.1, it follows that $(A,B)\in d^m(I)$. Again, since
$$
S^{m+n-2}=\left(\begin{array}{clcr}m+n-2\\n-1\end{array}\right)\triangle_{A,B}^{m-1}(I)B^{n-1}\otimes N^{n-1},
$$
$$
T^{m+n-2}=\left(\begin{array}{clcr}m+n-2\\n-1\end{array}\right)\delta^{m-1}_{A,B}(I)\otimes N^{n-1}
$$
and
$$
\triangle_{A,B}^{m-1}(I)B^{n-1}=0\Longrightarrow \triangle_{A,B}^{m-1}(I)=0,
$$
we must have $\triangle_{A,B}^{m-1}(I)\neq 0\neq \delta_{A,B}^{m-1}(I)$, i.e., $(A,B)\in$ strict-$d^m(I)$.
\end{demo}

\section {\sfstp Results: Commuting $A$ and $B$} Pairs of operators $(A,B)\in d^m(I)$ for which $[A,B]=0$ have a particular representation: There exists an $m$-nilpotent $N$ satisfying $[B,N]=0$ and such that $A=B^{-1}+N$ if $d=\triangle$ and $A=B+N$ if $d=\delta$. The following theorem is a generalization of (\cite{G}, [Prpoposition 8,9).

\begin{thm}\label{thm41} (a) Given $A, B\in \b$ such that $[A,B]=0$, if :

\vskip4pt\noindent (i) $(A,B)\in \triangle^m(I)$, then $B$ is invertible and there exists an $m$-nilpotent operator $N$ satisfying $[B,N]=0$ and $A=B^{-1}+N$;

\vskip4pt\noindent (ii) $(A,B)\in \delta^m(I)$, then there exists an $m$-nilpotent operator $N$ satisfying $[B,N]=0$ and $B=A+N$.

(b) If $A, B\in \b$ satisfy

\vskip4pt\noindent (iii) $\triangle_{A,B}^2(I)=0$, then $B$ is invertible if and only if $[A,B]=0$;

\vskip4pt\noindent (iv) $\delta_{A,B}^2(I)=0$, then $[A,B]=0$.
\end{thm}

\begin{demo} (i) The hypotheses $[A,B]=0$ and $(A,B)\in \triangle^m(I)$ imply $\triangle_{B,A}^m(I)=0=\triangle_{A,B}^m(I)$, and this in turn implies that both $A$ and $B$ are invertible. ($\triangle_{A,B}^m(I)=0$ implies $B$ is left $m$-invertible, hence left invertible, and $A$ is right $m$-invertible, hence right invertible). The invertibility of $B$  implies
$$
\triangle_{A,B}^m(I)=0 \Longleftrightarrow \sum_{j=0}^{m}(-1)^j\left(\begin{array}{clcr}m\\j\end{array}\right)A^jB^{-m+j}=0 \Longleftrightarrow \delta_{A,B^{-1}}^m(I)=0.
$$
Define $N\in \b$ by
$$
\delta_{A,B^{-1}}=L_A-R_{B^{-1}}=L_N.
$$
Then
$$
N^m=L_N^m(I)=\delta_{A,B^{-1}}^m(I)=0
$$
and
$$
L_N(I)=\delta_{A,B^{-1}}(I) \Longleftrightarrow A=B^{-1}+N.
$$
Evidently,
$$
I+NB=AB=BA=I+BN \Longleftrightarrow [B,N]=0.
$$

(ii) If $\delta_{A,B}^m(I)=0$, then define $N$ by $\delta_{A,B}=L_N$. We have
$$
A-B=\delta_{A,B}(I)=L_N(I)=N \Longleftrightarrow A=B+N,
$$
$$
N^m=L_N^m(I)=\delta_{A,B}^m(I)=0
$$
and
$$
NB=AB-B^2=BA-B^2=B(A-B)=BN.
$$

(iii) It being clear that if $[A,B]=0$ and $\triangle_{A,B}^2(I)=0$, then $B$ is invertible, we prove the converse. Assume $B$ is invertible. Then
$$
\triangle_{A,B}^2(I)=0 \Longleftrightarrow A^2-2AB^{-1}+B^{-2}=\delta_{A,B^{-1}}^2(I)=0.
$$
The operator $\delta_{A,B^{-1}}=L_N$ satisfies $\delta_{A,B^{-1}}^2(I)=L_N^2(I)=N^2=0$. Since $\delta_{A,B^{-1}}(I)=A-B^{-1}=N$, $(A-B^{-1})^2=A^2-AB^{-1}-B^{-1}A+B^{-2}=0$. Hence (since $\delta_{A,B^{-1}}^2(I)=A^2-2AB^{-1}+B^{-2}$)
$$
B^{-1}A=AB^{-1} \Longleftrightarrow [A,B]=0.
$$

(iv) If $\delta_{A,B}^2(I)=0$, then the operator $\delta_{A,B}=L_N$ satisfies $N=A-B$ and  $\delta_{A,B}^2(I)=N^2=0$. Hence
$$
\delta_{A,B}^2(I)=A^2-2AB+B^2=N^2=(A-B)^2=A^2-AB-BA+B^2 \Longrightarrow [A,B]=0.
$$
\end{demo}

If $\X={\mathbb{C}}^n$ is a finite dimensional Hilbert space, then $(A,B)\in\triangle^m(I)$ implies $B$ (also, $A$) is invertible. The conclusion of Theorem \ref{thm41} (iii) does not  extend to $m\geq 3$ even for $A,B\in B({\mathbb{C}}^n)$ (see \cite{G}, Example 10). Also, $[A,B]=0$ in (iv) of Theorem \ref{thm41} does not guarantee $\delta_{A,B}^2(I)=0$: Consider, for example, the operators $A$ and $DAD$ of Example \ref{ex22}: $[A,DAD]=0$ and $\delta_{A, DAD}^2(I) \neq 0$.

\

\noindent If  $B\in\B$  is  an algebraic operator, then $B$ has a representation $B=\oplus_{i=1}^r B|_{H_0(B-\lambda_i I)}$, where $\sigma(B)=\{\lambda_i:1\leq i\leq r\}$ and $H_0(B-\lambda_i I)=\{x\in\H:\lim_{n\rightarrow\infty}||(B-\lambda_iI)^nx||^{\frac{1}{n}}=0\}=(B-\lambda_i I)^{-{p_i}}(0)$ for some integer $p_i>0$.  ( Indeed, the points $\lambda_i$ are poles of the resolvent of $B$ of some order $p_i$.) Each $B_i=B|_{H_0(B-\lambda_i I)}$ has a representation $B_i=\lambda_iI_1+N_i$ for some nilpotent  $N_i$ (of order $p_i$) and $B=\oplus_{i=1}^r{\lambda_i I_i}+ \oplus_{i=1}^r{N_i}=B_0+N$ for some  normal operator $B_0$ and a nilpotent $N$ such that $[B_0,N]=0$. If we now assume (additionally) that $B$ is  $m$-isometric, then $B$ is invertible, hence $|\lambda_i|=1$ for all $\lambda_i\in\sigma(B)$ (\cite{FB}, Proposition 2.3). The operator $B_0$ being normal is (hence) unitary. Consequently, $(B^*_0,B_0)\in\triangle(I)$. Since $B$ is $m$-isometric (by hypothesis), it follows from an application of Corollary \ref{cor33} that the order $n$ of nilpotency of $N$ is given by $n=\frac{m+1}{2}$. Observe here that the invertibility of $B$ forces $m$ to be an odd integer (\cite{FB}, Proposition 2.4). We have proved:
\begin{pro}\label{pro41} Let  $B\in\B$ be an algebraic $m$-isometric operator. Then $m$ is odd and $B$ is the perturbation of a unitary operator by a commuting $\frac{m+1}{2}$-nilpotent operator.\end{pro}

Recall that operators $B\in B(\mathbb{C}^t)$, $\mathbb{C}^t$ the $t$-dimensional complex space, are algebraic. Proposition \ref{pro41} was proved for operators $B\in B(\mathbb{C}^t)$ by Berm\'{u}dez {\em et al} (\cite{BMN1}, Theorem 2.7). That an $m$-isometric operator $B\in B(\mathbb{C}^t)$ is the perturbation of a unitary by a commuting nilpotent was first observed by Agler {\em et al} (\cite{AHS}, Page 134).

%%%%%%%%%%%%%%%%%%%%%%%%%%%%%%%%%%%%%%%%%%%%%%%%%%%%%%%%%  REFERENCES

\vskip10pt \noindent\normalsize\rm B.P. Duggal, 8 Redwood Grove, London W5 4SZ, England (U.K.).\\
\noindent\normalsize \tt e-mail: bpduggal@yahoo.co.uk

\vskip6pt\noindent \noindent\normalsize\rm I. H. Kim, Department of
Mathematics, Incheon National University, Incheon, 22012, Korea.\\
\noindent\normalsize \tt e-mail: ihkim@inu.ac.kr

\end{document}